\documentclass[envcountsame]{llncs}

\usepackage{amsfonts}
\usepackage{amssymb}
\usepackage{mathtools}
\usepackage{tikz}
\usepackage{subfig}
\usepackage[ruled,vlined]{algorithm2e}
\usepackage{cite}

\newcommand{\seq}[1]{\langle #1 \rangle}
\newcommand{\intset}[2]{\{ #1 \;;\; #2 \}}
\newcommand{\extset}[1]{\{ #1 \}}
\newcommand{\cut}[2]{[\, #1 \mid #2 \,]}
\newcommand{\Goedel}{\mathfrak{g}}
\newcommand{\cantor}{2^\kappa}
\newcommand{\baire}{\kappa^\kappa}

\DeclareMathOperator{\Nat}{\ensuremath{\mathbb{N}}}
\DeclareMathOperator{\Rea}{\ensuremath{\mathbb{R}}}
\DeclareMathOperator{\Raz}{\ensuremath{\mathbb{Q}}}
\DeclareMathOperator{\No}{\ensuremath{No}}
\DeclareMathOperator{\Cut}{\ensuremath{Cut}}
\DeclareMathOperator{\dom}{\mathrm{dom}}

\title{Towards computable analysis on the\\ generalised real line}
\author{Lorenzo Galeotti\inst{1} \and Hugo Nobrega\inst{2}}

\institute{Fachbereich Mathematik, Universit\"at Hamburg,\\ Bundesstra{\ss}e 55, 20146 Hamburg, Germany \\
\email{lorenzo.galeotti@gmail.com}
\and
Institute for Logic, Language and Computation, Universiteit van Amsterdam, Postbus 94242, 1090 GE Amsterdam, The Netherlands\\
\email{h.nobrega@uva.nl}
}

\begin{document}
\maketitle
\begin{abstract}
	In this paper we use infinitary Turing machines with tapes of length $\kappa$ and which run for time \( \kappa \) as presented, e.g., by Koepke \& Seyfferth, to generalise the notion of type two computability to $\cantor$, where $\kappa$ is an uncountable cardinal with $\kappa^{<\kappa}=\kappa$.
	Then we start the study of the computational properties of $\Rea_\kappa$, a real closed field extension of $\Rea$ of cardinality $\cantor$, defined by the first author using surreal numbers and proposed as the candidate for generalising real analysis.
	In particular we introduce representations of $\Rea_\kappa$ under which the field operations are computable.
	Finally we show that this framework is suitable for generalising the classical Weihrauch hierarchy.
	In particular we start the study of the computational strength of the generalised version of the Intermediate Value Theorem. 
\end{abstract}

\section{Introduction}

The classical approach of computability theory is to define a notion of computability over $\omega$ and then extend that notion to any countable space via coding.
A similar approach is taken in computable analysis, where one usually defines a notion of computability over Cantor space $2^{\omega}$ or Baire space \( \omega^\omega \) by using the so-called type two Turing machines (T2TMs), and then extends that notion to spaces of cardinality at most the continuum via representations.
Intuitively a T2TM is a Turing machine in which a successful computation is one that runs forever (i.e., for $\omega$ steps). 
Using these machines one can compute functions over $2^{\omega}$, by stipulating that a function $f:2^{\omega}\rightarrow 2^{\omega}$ is computable if there is a T2TM which, when given $p\in \dom(f)$ as input, writes $f(p)$ on the output tape in the long run.
As an example, it is a classical result of computable analysis that, given the right representation of $\Rea$, the field operations are computable.
For an introduction to computable analysis we refer the reader to~\cite{Weihrauch}. 

Another classical application of T2TMs is the Weihrauch theory of reducibility (see, e.g.,~\cite{Brattka} for an introduction). 
The main aim of this theory is the study of the computational content of theorems of real analysis.
Since many of these theorems are of the form $\forall x \in \mathcal{X} \exists y \in \mathcal{Y} \, \varphi(x,y)$, with $\varphi(x,y)$ a quantifier free formula, they can be thought of as their own Skolem functions.
Given representations of \( \mathcal{X} \) and \( \mathcal{Y} \), Weihrauch reducibility provides a tool for comparing the computational strength of such functions, and therefore of the theorems themselves.
Using this framework, theorems from real analysis can be arranged in a complexity hierarchy analogous to the hierarchy of problems one has in classical computability theory.   

Recently, the study of the descriptive set theory of the generalised Baire spaces \( \baire \) and Cantor spaces $\cantor$ for cardinals $\kappa>\omega$ has been of great interest to set theorists. 
In~\cite{CIE2016} the second author provided the foundational basis for the study of \emph{generalised computable analysis}, namely the generalisation of computable analysis to generalised Baire and Cantor spaces.
In particular, in~\cite{CIE2016} the second author introduced $\Rea_\kappa$, a generalised version of the real line, and proved a version of the intermediate value theorem (IVT) for that space.

This paper is a continuation of~\cite{Thesis,CIE2016}, strengthening their results and answering in the positive the open question from~\cite{CIE2016} of whether a natural notion of computability exists for \( \cantor \).
We generalise the framework of type two computability to uncountable cardinals $\kappa$ such that $\kappa^{<\kappa}=\kappa$.
Then we use this framework to induce a notion of computability over the generalised real line $\Rea_\kappa$, showing that, as in the classical case, by using suitable representations, the field operations are computable.
Finally we will generalise Weihrauch reducibility to spaces of cardinality $\cantor$ and extend a classical result by showing that the generalised version of the IVT introduced in~\cite{CIE2016} is Weihrauch equivalent to a generalised version of the boundedness principle $\mathrm{B}_{\mathrm{I}}$. 

Throughout this paper $\kappa$ will be a fixed uncountable cardinal, as usual assumed to satisfy $\kappa^{<\kappa} = \kappa$, which in particular implies that $\kappa$ is a regular cardinal. 
The generalised Baire and Cantor spaces are equipped with their bounded topologies, i.e., the ones generated by the sets of the form \( \intset{x \in \lambda^\kappa}{\sigma \subset x } \) for \( \sigma \in \lambda^{<\kappa} \) and \( \lambda = 2 \) or \( \lambda = \kappa \), respectively.

\section{The surreal numbers}\label{section:Surreal}

The following definition as well as most of the results in this section are due to Conway~\cite{Conway} and have also been deeply studied by Gonshor in~\cite{Gonshor}.

A \emph{surreal number} is a function from an ordinal $\alpha$ to $\extset{+,-}$, i.e., a sequence of pluses and minuses of ordinal length.
We denote the class of surreal numbers by $\No$, and the set of surreal numbers of length strictly less than $\alpha$ by $\No_{<\alpha}$.
The \emph{length} of a surreal number $x$, denoted \( \ell(x) \), is its domain.
For surreal numbers $x$ and $y$, we define \( x < y \) if there exists \( \alpha \) such that \( x(\beta) = y(\beta) \) for all \( \beta < \alpha \), and (i) \( x(\alpha) = - \) and either \( \alpha = \ell(y) \) or \( y(\alpha) = + \), or (ii) \( \alpha = \ell(x) \) and \( y(\alpha) = + \).

In Conway's original idea, every surreal number is generated by filling some gap between shorter numbers.
The following theorem connects this intuition to the surreal numbers as we have defined them.
First, given sets of surreal numbers \( X \) and \( Y \), we write \( X < Y \) if for all \( x \in X \) and \( y \in Y \) we have \( x < y \).

\begin{theorem}[Simplicity theorem]\label{Theo:Simplicity}
	If $L$ and $R$ are two sets of surreal numbers such that $L<R$, then there is a unique surreal $x$ of minimal length such that $L < \extset{x} < R$, denoted by \( \cut{L}{R} \).
	Furthermore, for every $x\in \No$ we have \( x = \cut{L}{R} \) for \( L = \intset{y \in \No}{x>y\land y\subset x} \) and \( R = \intset{y \in \No}{x<y\land y\subset x} \).
	The pair \( \seq{L,R} \) is called the \emph{canonical cut} of \( x \).
\end{theorem}

Using the simplicity theorem Conway defined the field operations \( {+_\mathrm{s}} \), \( \cdot_\mathrm{s} \), \( -_\mathrm{s} \), and the multiplicative inverse over $\No$ and proved that these operations satisfy the axioms of real closed fields. 
These operations satisfy the following, where for any operation \( * \), surreal \( z \), and sets \( X, Y \) of surreals we use the notations \( z * X := \intset{z * x}{x \in X} \) and \( X * Y := \intset{x * y}{x \in X \text{ and } y \in Y} \).

\begin{theorem}\label{Theo:SurrealOperations}
	Let \( x = \cut{L_x}{R_x} \), \( y = \cut{L_y}{R_y} \) be surreal numbers.
	We have
	\[
	\begin{array}{rcl}
		x+_\mathrm{s}y
		&=& \cut{L_x+_\mathrm{s} y,x+_\mathrm{s} L_y}{R_x+_\mathrm{s} y,x+_\mathrm{s} R_y} \\
		\mathop{-_\mathrm{s}} x
		&=& \cut{\mathop{-_\mathrm{s}} R_x}{\mathop{-_\mathrm{s}} L_x} = \cut{\intset{\mathop{-_\mathrm{s}} x_R}{x_R\in R_x}}{\intset{\mathop{-_\mathrm{s}} x_L}{x_L\in L_x}} \\
		x\cdot_\mathrm{s} y
		&=& [L_x\cdot_\mathrm{s} y+_\mathrm{s} x\cdot_\mathrm{s} L_y -_\mathrm{s} L_x\cdot_\mathrm{s} L_y,R_x\cdot_\mathrm{s} y+_\mathrm{s} x\cdot_\mathrm{s} R_y  -_\mathrm{s} R_x\cdot_\mathrm{s} R_y \\
		&&\,{\mid}\, L_x\cdot_\mathrm{s} y+_\mathrm{s} x\cdot_\mathrm{s} R_y -_\mathrm{s} L_x\cdot_\mathrm{s}  R_y,R_x\cdot_\mathrm{s} y+_\mathrm{s} x\cdot_\mathrm{s} L_y -_\mathrm{s} R_x\cdot_\mathrm{s} L_y]
	\end{array}
	\]
	Now let \( z = \cut{L_z}{R_z} \) be a positive surreal number.
	Let \( r_{\seq{}} := 0 \) and recursively for every \( z_0,\ldots,z_n\in (L_z \cup R_z)\setminus \extset{0} \) let \( r_{\seq{z_0,\ldots, z_n}} \) be the solution for \( x \) of the equation $(z \mathop{-_\mathrm{s}} z_n) \cdot_\mathrm{s} r_{\seq{z_0,\ldots, z_{n-1}}} +_\mathrm{s} z_n \cdot_\mathrm{s} x = 1$.
	Then we have \( \frac{1}{z} = \cut{L'}{R'} \), where \( L' = \intset{ r_{\seq{z_0,\ldots, z_n}}}{n\in \Nat \) and \( z_i\in L_z \) for even-many \( i \leq n } \) and \( R' = \intset{ r_{\seq{z_0,\ldots, z_n}}}{n\in \Nat \) and \( z_i\in L_z \) for odd-many \( i \leq n } \).
\end{theorem}

On ordinals, the operations \( +_\mathrm{s} \) and \( \cdot_\mathrm{s} \) are the so-called \emph{natural} or \emph{Hessenberg} operations.
In particular, for any ordinal \( \alpha \) and natural number \( n \), we have \( \alpha +_\mathrm{s} n = \alpha + n \).

\section{The generalised real line}

A crucial property of the real line is its Dedekind completeness, forming the cornerstone of many theorems in real analysis.
However, it is a classical theorem that there are no real closed proper field extensions of $\Rea$ which are Dedekind complete (see, e.g., \cite[Theorem~8.7.3]{Cohn}).
We therefore need to replace Dedekind completeness with a weaker property. This was done in~\cite{Thesis,CIE2016}, and we repeat the central definitions here.

Let $X$ be an ordered set and $\kappa$ be a cardinal.
We say that $X$ is an \emph{$\eta_{\kappa}$-set} if whenever $L,R\subseteq X$ are such that $L<R$ and $\lvert L \cup R\rvert<\kappa$, there is $x\in X$ such that $L<\extset{x}<R$.
Let $K$ be an ordered field.
We call $\seq{L,R}$ a \emph{cut} over $K$ if $L,R\subseteq K$ and $L<R$.
Moreover we say that $\seq{L,R}$ is a \emph{Veronese cut} if it is a cut and $L$ has no maximum, $R$ has no minimum and for each $\varepsilon\in K^{+}$ there are $\ell\in L$ and $r\in R$ such that $r<\ell+\varepsilon$.
We say that $K$ is \emph{Veronese complete} if for each Veronese cut $\seq{L,R}$ there is $x\in K$ such that $L<\extset{x}<R$. 
Note that Veronese completeness is a reformulation of Cauchy completeness in terms of cuts (see, e.g.,~\cite{Ehrlich97}), so we can define the Cauchy completion of $\No_{<\kappa}$ as follows.

\begin{definition}
	$\Rea_\kappa=\No_{<\kappa}\cup \intset{\cut{L}{R}}{\seq{L, R} \text{ is a Veronese cut over $\No_{<\kappa}$}}.$
\end{definition}
\begin{theorem}[Galeotti~\cite{CIE2016}]\label{Theo:RkProp}
	The field $\Rea_\kappa$ is the unique Cauchy-complete real closed field extension of $\Rea$ which is an \( \eta_\kappa\)-set of cardinality $2^\kappa$, degree \( \kappa \), and in which \( \No_{<\kappa} \) can be densely embedded.
\end{theorem}

In view of the previous theorem from now on we will call $\No_{<\kappa}$ the $\kappa$-rational numbers and we use the symbol $\Raz_\kappa$ instead of $\No_{<\kappa}$.


The field \( \Rea_\kappa \) is a suitable setting for generalising results from classical analysis.
For example, a generalised version of the intermediate value theorem~\cite{CIE2016}, a generalised version of the extreme value theorem~\cite{Thesis}, and recently a generalised version of the Bolzano-Weierstra{\ss} theorem (for \( \kappa \) weakly compact)~\cite{BonnSeminar2016} have been proved to hold for \( \Rea_\kappa \).
In this section we briefly recall some of the definitions from~\cite{CIE2016} which will be needed in the last part of this paper.


A \emph{$\kappa$-topology} over a set $X$ is a collection of subsets $\tau$ of $X$ satisfying: $\emptyset,X\in \tau$; for any \( \alpha < \kappa \), if \( \extset{A_i}_{i\in \alpha} \) is a collection of sets in $\tau$ then $\bigcup_{i<\alpha} A_i \in \tau$; and for all \( A,B \in \tau \), we have \( A\cap B\in \tau \).
With \( \kappa \)-topologies one can define direct analogues of many topological notions.
We refer to these with the prefix ``\( \kappa \)-''; thus we have \( \kappa \)-open sets, \( \kappa \)-continuous functions, \( \kappa \)-topologies generated by families of subsets of a set, etc.
Note that, unlike the classical case of the interval topology over \( \Rea \), the interval \( \kappa \)-topologies over \( \Rea_\kappa \) in which the intervals have endpoints in \( \Rea_\kappa \cup \extset{-\infty,+\infty} \) or in \( \Raz_\kappa \cup \extset{-\infty,+\infty} \) are different in general.
In what follows we will only consider the generalised real line \( \Rea_\kappa \) equipped with the former.


\begin{theorem}[$\mathrm{IVT}_\kappa$~\cite{CIE2016}]\label{Theo:IVT}
	Let $a,b \in \Rea_\kappa$ and $f:[0,1]\rightarrow \Rea_\kappa$ be a $\kappa$-continuous function.
	Then for every $r \in [f(0), f(1)]$ there exists $c\in [0,1]$ such that $f(c)=r$.
\end{theorem} 

\section{Generalised type two Turing machines}

In this section we define a generalised version of type two Turing machines (T2TMs). 
We will only sketch the definition of \( \kappa \)-Turing machines, which were developed by several people (e.g.,~\cite{Dawson,Koepke09,Rin}); we are going to follow the definition of Koepke and Seyfferth~\cite[\S\ 2]{Koepke09}.

A \emph{$\kappa$-Turing machine} has the following tapes of length $\kappa$: finitely many read-only tapes for the input, finitely many read and write scratch tapes and one write-only tape for the output.
Each cell of each tape has either \( 0 \) or \( 1 \) written in it at any given time, with the default value being \( 0 \).
These machines can run for infinite time of ordinal type $\kappa$; at successor stages of a computation a $\kappa$-Turing machine behaves exactly like a classical Turing Machine, while at limit stages the contents of each cell of each tape and the positions of the heads is computed using inferior limits. 

As in the classical case \( \kappa=\omega \), the difference between \( \kappa \)-Turing machines and type 2 \( \kappa \)-Turing machines is \emph{not} on the machinery level, but rather on the notion of what it means for a machine to compute a function.
A partial function \( f:2^{<\kappa}\to2^{<\kappa} \) is \emph{computed} by a \( \kappa \)-Turing machine \( M \) if whenever \( M \) is given \( x \in \dom(f) \) as input, its computation halts after fewer than \( \kappa \) steps with \( f(x) \) written on the output tape.
A partial function \( f:\cantor\to\cantor \) is \emph{type two-computed} by a \( \kappa \)-Turing machine \( M \), or \emph{computed by the type 2 \( \kappa \)-Turing machine \( M \)}, or simply \emph{computed} by \( M \), if whenever \( M \) is given \( x \in \dom(f) \) as input, for every \( \alpha < \kappa \) there exists a stage \( \beta < \kappa \) of the computation at which \( f(x) \upharpoonright \alpha \) is written on the output tape.
We abbreviate \emph{type 2 \( \kappa \)-Turing machine} by \emph{T2\( \kappa \)TM}.
%
An \emph{oracle T2$\kappa$TM} is a T2$\kappa$TM with an additional read-only input tape of length \( \kappa \), called its \emph{oracle tape}.
A partial function \( f: \cantor \to \cantor \) is \emph{computable with an oracle} if there exists an oracle T2\(\kappa\)TM \( M \) and \( x \in \cantor \) such that \( M \) computes \( f \) when \( x \) is written on the oracle tape.
Note that by minor modifications of classical proofs one can prove that T2$\kappa$TMs are closed under recursion and composition, and that there is a universal T2\( \kappa \)TM. 
In what follows, the term \emph{computable} will mean \emph{computable by a T2\(\kappa\)TM}, unless specified otherwise.


%


\begin{theorem}
	A partial function \( f:\cantor \to \cantor \) is continuous iff it is computable with some oracle.
	\label{Theo:Oracle}
\end{theorem}

\section{Represented spaces}

In this section we generalise the classical definitions of the theory of represented spaces to \( \cantor \) (see, e.g.,~\cite{Weihrauch,Arno} for the classical case).

A \emph{represented space} $\mathbf{X}$ is a pair $(X,\delta_X)$ where $X$ is a set and $\delta_X:\cantor\rightarrow X$ is a partial surjective function.
As usual a multi-valued function between represented spaces is a multi-valued function between the underlying sets. 
Let $f:\mathbf{X}\rightrightarrows \mathbf{Y}$ be a partial multi-valued function between represented spaces.
We call $F: \cantor \rightarrow \cantor$ a \emph{realizer} of $f$, in symbols $F\vdash f$, if for every $x\in \dom(\delta_X)$ we have that $\delta_Y(F(x))\in f(\delta_X(x))$.
Given a class \( \Gamma \) of functions between \( \cantor \) and \( \cantor \), we say $f$ is $(\delta_X,\delta_Y) \text{-} \Gamma$, or \( \delta_X \text{-} \Gamma \) in case \( \delta_X = \delta_Y \), if \( f \) has a realizer in \( \Gamma \). 
For example, a function $f:X\rightarrow Y$ is $(\delta_X,\delta_Y) \text{-}$computable if it has a computable realizer. 

Let $f$ and $g$ be two multi-valued functions between represented spaces.
Then we say that $f$ is \emph{strongly topologically-Weihrauch reducible} to $g$, in symbols $f\leq^\mathrm{t}_\mathrm{W} g$, if there are two continuous functions $H,K:\cantor\rightarrow \cantor$ such that $H\circ G\circ K\vdash f$ whenever $G\vdash g$.
If the functions \( H,K \) above can be taken to be computable, then we say \( f \) is \emph{strongly Weihrauch reducible} to $g$, in symbols $f\leq_\mathrm{W} g$.%
\footnote{Carl has also introduced a notion of generalized (strong) Weihrauch reducibility in~\cite{Carl}.
Because his goal is to investigate multi-valued (class) functions on \( V \), the space of codes he uses is the class of ordinal numbers, considered with the ordinal Turing machines of Koepke~\cite{Koepke}.
Therefore his approach is significantly different from ours, and we do not know of any connections between the two.}
As usual, if $f\leq^\mathrm{t}_\mathrm{W} g$ and $g\leq^\mathrm{t}_\mathrm{W} f$ then we say that $f$ is \emph{strongly topologically-Weihrauch equivalent} to $g$ and write $f\equiv^\mathrm{t}_{\mathrm{W}} g$.
The relation \( \equiv_{\mathrm{W}} \) is defined analogously.


Let $\delta: \cantor\rightarrow X$ and $\delta':\cantor\rightarrow X$ be two representations of a space $X$.
Then we say that $\delta$ \emph{continuously reduces} to $\delta'$, in symbols $\delta\leq_\mathrm{t} \delta'$, if there is a continuous function $h:\cantor\rightarrow \cantor$ such that for every $p\in \dom(\delta)$ we have $\delta(p)=\delta'(h(p))$.
Similarly we say that $\delta$ \emph{computably reduces} to $\delta'$, in symbols $\delta\leq\delta'$, if \( h \) above can be taken computable.
If $\delta\leq_\mathrm{t}\delta'$ and $\delta'\leq_\mathrm{t} \delta$ we say that $\delta$ and $\delta'$ are continuously equivalent and write $\delta\equiv_\mathrm{t} \delta'$, and similarly for the computable case.
Note that as in classical computable analysis if $\delta \leq \delta'$ and $f$  is $\delta$-computable then $f$ is also $\delta'$-computable.  
Finally, as in the classical case, given two represented spaces $\mathbf{X}$ and $\mathbf{Y}$, we can define canonical representations for the product space $\mathbf{X}\times \mathbf{Y}$, the union space $\mathbf{X}+\mathbf{Y}$ and the space of continuous functions $[\mathbf{X}\rightarrow\mathbf{Y}]$.
In particular, as in classical computable analysis $[\mathbf{X}\rightarrow\mathbf{Y}]$ can be represented as follows:
$\delta_{[X \to Y]}(p)=f$
iff $p=0^n1 p'$ with $p'\in \cantor$ and $n\in \Nat$ is a code for an oracle T2$\kappa$TM which $(\delta_X,\delta_Y)$-computes $f$ when given the oracle $p'$. 

Recall that the following relation is a well-ordering of the class of pairs of ordinal numbers: \( \seq{\alpha_0,\beta_0} \prec \seq{\alpha_1,\beta_1} \) iff \( \seq{\max (\alpha_0,\beta_0),\alpha_0,\beta_0} \) is lex\-i\-co\-graph\-i\-cal\-ly\--less than \( \seq{\max (\alpha_1,\beta_1),\alpha_1,\beta_1} \).
The \emph{G\"odel pairing function} is given by \( \Goedel(\alpha,\beta) = \gamma \) iff \( \seq{\alpha,\beta} \) is the \( \gamma \)\textsuperscript{th} element in \( \prec \).
Given sequences \( \seq{w_\alpha}_{\alpha<\kappa} \) and \( \seq{p_\alpha}_{\alpha<\beta} \) of elements in \( 2^{<\kappa} \) and \( \cantor \), respectively, we define elements \( q := [w_\alpha]_{\alpha<\kappa} \) and \( p := (p_\alpha)_{\alpha<\kappa} \) in \( \cantor \) by letting \( q \) be the concatenation of the \( w_\alpha \) and \( p(\Goedel(\alpha,\beta)) = p_\alpha(\beta) \).

We fix the following representations of $\kappa$ and \( \baire \): $\delta_\kappa(p)=\alpha$ iff $p=0^{\alpha}1\mathbf{0}$, where $\mathbf{0}$ is the constant $0$ $\kappa$-sequence, $\delta_{\baire}(p) = x$ iff $p = [0^{\alpha_\beta+1}1]_{\beta < \kappa} \text{ and } x = \seq{\alpha_\beta}_{\beta < \kappa}$.
It is straightforward to see that a function \( f:\kappa\to\kappa \) is \( \delta_\kappa \)-computable iff it is computable by a \( \kappa \)-machine as in~\cite[Definition 2]{Koepke09}.

\begin{lemma}
	The restriction of \( \Goedel \) to \( \kappa \times \kappa \) is a \( \delta_\kappa \)-computable bijection between \( \kappa \times \kappa \) and \( \kappa \), and has a \( \delta_\kappa \)-computable inverse.
	\label{Goedel_pairing}
\end{lemma}

\begin{proposition}\label{Prop:DeltaKappaMinimal}
	$\delta_{\baire}$ is \( \leq \)-maximal among the continuous representations of~$\kappa^{\kappa}$.
\end{proposition}

\section{Representing $\Rea_\kappa$}

In classical computable analysis one can show that many of the natural representations of $\Rea$ are well behaved with respect to type two computability.
In this section we show that some of these results naturally extend to the uncountable case.
First we introduce representations for generalised rational numbers, which will serve as a starting point to representing $\Rea_\kappa$.
As we have seen in the introduction, surreal numbers can be expressed as binary sequences and, because of the simplicity theorem, as cuts.
It is then natural to introduce two representations which reflect this fact.
Let $p\in \cantor$ and $q\in \Raz_\kappa$.
We define $\delta_{\Raz_\kappa}(p)=q$ iff  $p=[w_\alpha]_{\alpha < \kappa}$ where $w_\alpha:=00$ if $\alpha\in \dom(q)$ and $q(\alpha)=-$, $w_\alpha:=01$ if $\alpha\notin \dom(q)$, and finally $w_\alpha:=11$ if $\alpha\in \dom(q)$ and $q(\alpha)=+$.
It is not hard to see that since every rational is a sequence of $+$ and $-$ of length less than $\kappa$ the function $\delta_{\Raz_\kappa}$ is indeed a representation of $\Raz_\kappa$.
Now we define a representation based on cuts by recursion on the simplicity structure of the surreal numbers. 
We define $\delta^0_{\Cut_{\Raz_\kappa}}(p)=0$ iff $p = (p_\alpha)_{\alpha<\kappa}$ and $p_\alpha=[10]_{\beta<\kappa}$ for every $\alpha<\kappa$.
For $\alpha>0$ we define $\delta^\alpha_{\Cut_{\Raz_\kappa}}(p)=\cut{L}{R}$ where $p = (p_\alpha)_{\alpha<\kappa}$ and:
\begin{enumerate}
	\item
		$p_{\alpha}\in \dom(\bigcup_{\gamma<\alpha}\delta^{\gamma}_{\Raz_\kappa})\cup \extset{[10]_{\beta < \kappa}}$ for every $\alpha < \kappa$,
	\item
		for all even\footnote{We call an ordinal \( \alpha \) \emph{even} if \( \alpha = \lambda + 2n \) for some limit \( \lambda \) and natural \( n \), \emph{odd} otherwise.} $\alpha<\kappa$, if $p_{\alpha}=[10]_{\beta < \kappa}$ then for all even $\beta>\alpha$ we have $p_{\beta}=[10]_{\beta < \kappa}$,
	\item
		for all odd $\alpha<\kappa$, if $p_{\alpha}=[10]_{\beta < \kappa}$ then for all odd $\beta>\alpha$ we have $p_{\beta}=[10]_{\beta < \kappa}$,
	\item finally: $L=\intset{\delta^{\gamma}_{\Cut_{\Raz_\kappa}}(p_\beta)}{\text{ $\gamma<\alpha$, $\beta<\kappa$ is even and $p_\beta\in \dom(\delta^\gamma_{\Cut_{\Raz_\kappa}})$} }$ and 
		$R=\intset{\delta^{\gamma}_{\Cut_{\Raz_\kappa}}(p_\beta)}{\text{$\gamma<\alpha$, $\beta<\kappa$ is odd and $p_\beta\in \dom(\delta^\gamma_{\Cut_{\Raz_\kappa}})$} }$.
\end{enumerate}
Then we define $\delta_{\Cut_{\Raz_\kappa}}:=\bigcup_{\gamma<\kappa}\delta^\gamma_{\Cut_{\Raz_\kappa}}$. 

Note that $\delta_{\Cut_{\Raz_\kappa}}$ is surjective, since for every \( x \in \Raz_\kappa \) there exists \( p \in \dom(\delta_{\Cut_{\Raz_\kappa}}) \) such that \( \delta_{\Cut_{\Raz_\kappa}}(p) \) is the canonical cut for \( x \). 
Therefore \( \delta_{\Cut_{\Raz_\kappa}} \) is indeed a representation of \( \Raz_\kappa \).


\begin{lemma}\label{Lemma:RatRepEq}
	$\delta_{\Raz_\kappa}\equiv \delta_{\Cut_{\Raz_\kappa}}$.
\end{lemma}
\begin{proof}
	First we show that $\delta_{\Raz_\kappa}\leq \delta_{\Cut_{\Raz_\kappa}}$.
	Let $p\in \dom(\delta_{\Raz_\kappa})$.
	The conversion can be done recursively.
	If $p$ is a code for the empty sequence\footnote{Note that this can be checked just by looking at the first two bits of $p$.} we just return a representation for $\cut{\emptyset}{\emptyset}$.
	Otherwise we compute two subsets $L_s:= \intset{p'01}{p'11\subset p}$ and $R_s:=\intset{p'01}{p'00\subset p}$.
	Then we compute recursively the cuts for the elements of $L_s$ and $R_s$ and return them respectively as the left and right sets of the cut representation of $p$.
	It easy to see that the algorithm computes a code for the canonical cut of $\delta_{\Raz_\kappa}(p)$. 

	Now we will show that $\delta_{\Cut_{\Raz_\kappa}} \leq \delta_{\Raz_\kappa}$.
	Let $p\in \dom(\delta_{\Cut_{\Raz_\kappa}})$.
	If $p$ is a code for the cut $\cut{\emptyset}{\emptyset}$ we return a representation of the empty sequence.
	If $p$ is the code for the cut $\cut{L}{R}\neq \cut{\emptyset}{\emptyset} $.
	We first recursively compute the sequences for the element of $L$ and $R$, call the sets of these sequences $L_s$ and $R_s$.
	Now suppose $\alpha<\kappa$ is even and we want to compute the value at $\alpha$ and $\alpha+1$ of the output sequence.
	We first compute $M_L$ and $m_R$ respectively the minimal and maximal in $\extset{00,01,11}$ such that for every $p'\in L_s$ and $p''\in R_s$ we have $p'(\alpha)p'(\alpha+1)\leq M_L$ and $m_R\leq p''(\alpha)p''(\alpha+1)$.
	Then by a case distinction on $M_L$ and $m_R$ we can decide the $i$\textsuperscript{th} sign of the output.
	For example if the output is already smaller than $R_s$, $M_L=00$ (i.e. $-$) and $m_R=00$ (i.e. $-$) then we can output the sequence $01$ (i.e. undefined). All the other combinations can be treated similarly.
\end{proof}

\begin{lemma}
	The operations $+_\mathrm{s}$, $-_\mathrm{s}$, $\cdot_\mathrm{s}$, $\frac{1}{x}$ and the order $<$ are $\delta_{\Cut_{\Raz_\kappa}}$-computable.
\end{lemma}
\begin{proof}
	We will only prove the lemma for $+_\mathrm{s}$.
	Given $q,q'\in \Raz_\kappa$ we want to $\delta_{\Cut_{\Raz_\kappa}}$-compute $q+_\mathrm{s}q'$.
	The algorithm is given by recursion.
	If $q=0$ (similarly for $q'=0$)\footnote{Note that this is easily computable, it is in fact enough to check that $L$ and $R$ are empty, and this can be done just by checking the first two bits of the first sequence in the left and in the first sequence on the right.} copy the code of $q'$ on the output tape.
	If neither $q$ nor $q'$ are $0$ then by using Theorem \ref{Theo:SurrealOperations} we compute a representation for $q+_{\mathrm{s}}q'$ (note that this involves the computation of less than $\kappa$ many rational sums of shorter length). 
	Finally, since the resulting code would not in general be in $\dom(\delta_{\Cut_{\Raz_\kappa}})$, we use the algorithms of the previous lemma to convert $q+_\mathrm{s}q'$ to a sign sequence code and than we convert it back to an element in $\dom(\delta_{\Cut_{\Raz_\kappa}})$. 
	By using the second algorithm from the previous proof we can convert every element in $ L_{q+_\mathrm{s}q'}$ and  in $R_{q+_\mathrm{s}q'}$ into a sequence (note that by induction the codes of these cuts are in $\dom(\delta_{\Cut_{\Raz_\kappa}})$ so we can use the algorithm). 
	Then by the same method used in the previous lemma, we can compute the code of the sequence representation for $q+_\mathrm{s}q'$. 
	Once we have the code of the sequence representation for $q+_\mathrm{s}q'$ we can convert it to a code of the cut representation by using the first algorithm from the previous lemma.
\end{proof}

Given that \( \Rea_\kappa \) is the Cauchy completion of \( \Raz_\kappa \), the following is a natural representation of \( \Rea_\kappa \).
We let \( \delta_{\Rea_\kappa}(p) = x \) iff \( p = (p_\alpha)_{\alpha < \kappa} \), where for each \( \alpha < \kappa \) we have \( p_\alpha \in \dom(\delta_{\Raz_\kappa})\), \( \delta_{\Raz_\kappa}(p_{\alpha}) < x +_\mathrm{s} \frac{1}{\alpha+1} \), and \( x < \delta_{\Raz_\kappa}(p_{\alpha}) +_\mathrm{s} \frac{1}{\alpha+1} \).
It is routine to check the following.

\begin{theorem}
	The field operations $+_\mathrm{s}$, $-_\mathrm{s}$, $\cdot_\mathrm{s}$, and $\frac{1}{x}$ are $\delta_{\Rea_\kappa}$-computable.
	\label{R_kappa_operations_computable}
\end{theorem}
\begin{proof}
	Let us do the proof for \( \cdot_\mathrm{s} \), the others being similar.
	Given codes \( p = (p_\alpha)_{\alpha < \kappa} \) and \( q = (q_\alpha)_{\alpha < \kappa} \) for \( x, y \in \Rea_\kappa \) respectively, let \( x_\alpha = \delta_{\Raz_\kappa}(p_\alpha) \) and \( y_\alpha = \delta_{\Raz_\kappa}(q_\alpha) \).
	Note that for each \( \alpha \) we can compute some \( \alpha' \) such that \( \frac{1}{\alpha'+1}(x_0 +_\mathrm{s} y_0 +_\mathrm{s} 3) \leq \frac{1}{\alpha+1} \).
	We then output \( r = (r_\alpha)_{\alpha < \kappa} \), where \( r_\alpha \) is a \( \delta_{\Raz_\kappa} \)-name for \( x_{\alpha'}y_{\alpha'} \).

	We have $xy -_\mathrm{s} x_{\alpha'}y_{\alpha'}= x(y -_\mathrm{s} y_{\alpha'}) +_\mathrm{s} y_{\alpha'}(x -_\mathrm{s} x_{\alpha'})
	< (x_0 +_\mathrm{s} 1)\frac{1}{\alpha'+ 1} +_\mathrm{s} (y_{0}+_\mathrm{s}2)\frac{1}{\alpha'+ 1} 
	\leq \frac{1}{\alpha+ 1}$,
	as desired, and likewise we can prove \( x_{\alpha'}y_{\alpha'} -_\mathrm{s} xy < \frac{1}{\alpha+ 1} \).
\end{proof}

On the other hand, the following is suggested by the definition of \( \Rea_\kappa \) as the collection of Veronese cuts over \( \Raz_\kappa \).
We let \( \delta^{\mathrm{V}}_{\Rea_{\kappa}}(p) = x \) iff \( p = (p_\alpha)_{\alpha < \kappa} \), where for each \( \alpha < \kappa \) we have \( p_\alpha \in \dom(\delta_{\Raz_\kappa})\) and \( x = \cut{L}{R} \), with \( L = \{ \delta_{\Raz_\kappa}(p_\alpha) \;;\; \alpha < \kappa \text{ is even} \} \);	 \( R = \{ \delta_{\Raz_\kappa}(p_\alpha) \;;\; \alpha < \kappa \text{ is odd} \} \); and for each even \( \alpha < \kappa \) we have \( \delta_{\Raz_\kappa}(p_{\alpha+1}) < \delta_{\Raz_\kappa}(p_\alpha) +_\mathrm{s} \frac{1}{\alpha+ 1} \).
\begin{theorem}\label{Theo:RedRkappaRep}
	\( \delta_{\Rea_{\kappa}} \equiv \delta^{\mathrm{V}}_{\Rea_{\kappa}} \).
\end{theorem}
\begin{proof}
	To reduce \( \delta^{\mathrm{V}}_{\Rea_{\kappa}} \) to \( \delta_{\Rea_{\kappa}} \), given \( p = (p_\alpha)_{\alpha < \kappa} \), we output \( q = (q_\alpha)_{\alpha < \kappa} \) by making \( q_\alpha \) equal to \( p_\beta \), where \( \beta \) is the \( \alpha \)\textsuperscript{th} even ordinal.
	It is now easy to see that \( q \) is a \( \delta_{\Rea_{\kappa}} \)-name for \( \delta^{\mathrm{V}}_{\Rea_{\kappa}}(p) \).

	For the converse reduction, given \( p = (p_\alpha)_{\alpha < \kappa} \), we output \( q = (q_\alpha)_{\alpha < \kappa} \) where for each even \( \alpha \) we let \( q_\alpha \) be a \( \delta_{\Raz_\kappa} \)-name for \( \delta_{\Raz_\kappa}(p_{2 \cdot_{\mathrm{s}} \alpha+2}) -_\mathrm{s} \frac{1}{2 \cdot_{\mathrm{s}} \alpha+3} \) and \( q_{\alpha+1} \) be a \( \delta_{\Raz_\kappa} \)-name for \( \delta_{\Raz_\kappa}(p_{2 \cdot_{\mathrm{s}} \alpha+2}) +_\mathrm{s} \frac{1}{2 \cdot_{\mathrm{s}} \alpha+3} \).
	Then letting \( L := \{ \delta_{\Raz_\kappa}(p_\alpha) \;;\; \alpha < \kappa \text{ is even} \} \) and \( R := \{ \delta_{\Raz_\kappa}(p_\alpha) \;;\; \alpha < \kappa \text{ is odd} \} \) we have \( L < \extset{x} <  R \) and for each even \( \alpha < \kappa \)
	we have
	$\delta_{\Raz_\kappa}(q_{\alpha+1})
	= \delta_{\Raz_\kappa}(p_{2 \cdot_{\mathrm{s}} \alpha+2}) +_\mathrm{s} \frac{1}{2 \cdot_{\mathrm{s}} \alpha+3} 
	= \delta_{\Raz_\kappa}(q_{\alpha}) +_\mathrm{s} \frac{2}{2 \cdot_{\mathrm{s}} \alpha+3} 
	< \delta_{\Raz_\kappa}(q_{\alpha}) +_\mathrm{s} \frac{1}{\alpha+1}$,
	as desired.
\end{proof}



\section{Generalised boundedness principles and the IVT}

As shown in, e.g.,~\cite{Brattka,ArnoLowBasis}, the so-called \emph{boundedness principles} and \emph{choice principles} are important building blocks in characterizing the Weihrauch degrees of interest in computable analysis. 
In this section we focus on the study of IVT and its relationship with the boundedness principle $\mathrm{B}_{\mathrm{I}}$.
In particular we generalise a classical result from Brattka and Gherardi~\cite{Brattka}, proving that $\mathrm{IVT}_\kappa$ is Weihrauch equivalent to a generalised version of $\mathrm{B}_{\mathrm{I}}$.
This strengthens a result from~\cite{CIE2016}, namely that $\mathrm{B}_{\mathrm{I}}$ is continuously reducible to $\mathrm{IVT}_\kappa$.

The theorem $\mathrm{IVT}_\kappa$ as stated in Theorem~\ref{Theo:IVT} can be considered as the partial multi-valued function $\mathrm{IVT}_\kappa:\mathrm{C}_{[0,1]}\rightrightarrows[0,1]$ defined as follows: $\mathrm{IVT}_\kappa(f)=\intset{c\in [0,1]}{f(c)=0 }$, where $[0,1]$ is represented by $\delta_{\Rea_\kappa}{\upharpoonright} [0,1]$ and $\mathrm{C}_{[0,1]}$ is endowed with the standard representation of $[[0,1]\rightarrow \Rea_\kappa]$ restricted to $\mathrm{C}_{[0,1]}$.
By lifting the classical proof to $\kappa$ it is easy to show that this version of $\mathrm{IVT}_\kappa$ is not
continuous, and thus also not computable, relative to these representations.

To introduce the boundedness principle \( \mathrm{B}^{\kappa}_\mathrm{I} \), we will need the following represented spaces.
Let \( \mathbf{S}_{\mathrm{b}}^{\uparrow} \) be the space of bounded increasing sequences of \( \kappa \)-rationals, represented by letting \( p \) be a name for \( \seq{x_\alpha}_{\alpha < \kappa} \) iff \( p = (p_\alpha)_{\alpha<\kappa} \) where \( p_\alpha \in \dom_{\Raz_\kappa} \) and \( \delta_{\Raz_\kappa}(p_\alpha) = x_\alpha \) for each \( \alpha < \kappa \).
The represented space \( \mathbf{S}_{\mathrm{b}}^{\downarrow} \) is defined analogously, with bounded decreasing sequences of \( \kappa \)-rationals.
Note that, unlike the classical case of the real line, not all limits of bounded monotone sequences of length \( \kappa \) exist in \( \Rea_\kappa \).
Therefore, although for the real line the spaces \( \mathbf{S}_{\mathrm{b}}^{\uparrow} \) and \( \mathbf{S}_{\mathrm{b}}^{\downarrow} \) naturally correspond to the spaces of \emph{lower reals} \( \Rea_< \) and \emph{upper reals} \( \Rea_> \), respectively, in our generalised setting the correspondence fails. 
We define $\mathrm{B}^\kappa_\mathrm{I}$ as the principle which, given an increasing sequence $\seq{q_\alpha}_{\alpha < \kappa}$ and decreasing sequence $\seq{q'_\alpha}_{\alpha < \kappa}$ in $\Raz_\kappa$ for which there exists  $x\in \Rea_\kappa$ such that 
$\intset{q_\alpha}{\alpha < \kappa} \leq \extset{x} \leq \intset{q'_\alpha}{\alpha < \kappa}$, picks one such \( x \).
Formally we have the partial multi-valued function \( \mathrm{B}^{\kappa}_\mathrm{I}:\mathbf{S}_{\mathrm{b}}^{\uparrow} \times \mathbf{S}_{\mathrm{b}}^{\downarrow} \rightrightarrows \Rea_\kappa \) with \( x \in \mathrm{B}^{\kappa}_\mathrm{I}(s, s') \) iff \( \intset{s(\alpha)}{\alpha < \kappa} \leq \extset{x} \leq \intset{s'(\alpha)}{\alpha < \kappa}. \)

\begin{lemma}
	Let \( f: [0,1] \to \Rea_\kappa \) and \( x \in \Rea_\kappa \).
	Suppose there exists a sequence \( \seq{x_\alpha}_{\alpha < \kappa} \) of pairwise distinct elements of \( [0,1] \) such that \( f(x_\alpha) = x \) if \( \alpha < \kappa \) is even and \( f(x_\alpha) \neq x \) otherwise, and such that for any odd \( \alpha, \beta < \kappa \) there exists an even \( \gamma < \kappa \) such that \( x_\gamma \) is between \( x_\alpha \) and \( x_\beta \).
	Then \( f \) is not \( \kappa \)-continuous.
	\label{kappa_continuous_cannot_flip_kappa_times}
\end{lemma}
\begin{proof}
	If such a sequence exists, then either the preimage of the \( \kappa \)-open set \( (x,+\infty) \) or of the \( \kappa \)-open set \( (-\infty,x) \) under \( f \) must contain \( x_\alpha \) for \( \kappa \)-many of the odd \( \alpha < \kappa \), and thus cannot be \( \kappa \)-open.
\end{proof}

\begin{lemma}\label{Theo:continuity}
	Let $f:[0,1]\rightarrow\Rea_\kappa$ be $\kappa$-continuous an let $\beta,\beta'<\kappa$, $y\in \Rea_\kappa$ and let $\seq{r_\alpha}_{\alpha < \beta}$ and $\seq{r'_\alpha}_{\alpha < \beta'}$ be two sequences in \( [0,1]\) such that $\intset{r_\alpha}{\alpha<\beta }<\intset{r'_\alpha}{\alpha<\beta'}$ and $\intset{f(r_\alpha)}{\alpha < \beta} < \extset{y} < \intset{f(r'_{\alpha})}{\alpha < \beta'}$.
	Then there is $x\in [0,1]$ such that $\intset{r_\alpha}{\alpha<\beta} < \extset{x} < \intset{r'_\alpha}{\alpha<\beta'}$ and $f(x)=y$.
\end{lemma}
\begin{proof}
	Assume not. 
	Without loss of generality we can assume that for every $x$ such that $\intset{r_\alpha}{\alpha<\beta}< \extset{x} <\intset{r'_\alpha}{\alpha<\beta'}$ we have $f(x)>y$ (a similar proof works for $f(x)<y$).
	Note that the set $\intset{r_\alpha}{\alpha<\beta}$ has cofinality at most $\beta<\kappa$ and, since $\Rea_\kappa$ is an $\eta_\kappa$-set, it follows that $R=\intset{r\in [0,1]}{\forall \alpha<\beta .\, r_\alpha<r}$ has coinitiality $\kappa$. 
	Therefore $R$ is not $\kappa$-open. Now since $f$ is $\kappa$-continuous we have that $f^{-1}[(y,+\infty)]$ is $\kappa$-open. 
	Therefore $f^{-1}[(y,+\infty)]=\bigcup_{\alpha\in \gamma}(y_\alpha,b_\alpha)$ with $\gamma<\kappa$ and $y_\alpha,b_\alpha\in 
	[0,1] $
	for every $\alpha<\gamma$.
	Now consider the set $I : =\intset{\alpha\in\gamma}{(y_\alpha,b_\alpha)\cap R\neq \emptyset}.$
	We have that $R\subset\bigcup_{\alpha\in I}(y_\alpha,b_\alpha)$. 
	Note that since $R$ is not $\kappa$-open we have $R\neq \bigcup_{\alpha\in I}(y_\alpha,b_\alpha)$. 
	Now assume $r\in \bigcup_{\alpha\in I}(y_\alpha,b_\alpha)\setminus R$, so that there is $\alpha\in I$ such that $r\in (y_\alpha,b_\alpha)$. 
	Take $r'\in (y_\alpha,b_\alpha)\cap R$. 
	By the fact that $r\notin R$, there is $\alpha'<\beta$ such that $r<r_{\alpha'}$ and by $\mathrm{IVT}_\kappa$ there is a root of $f$ between $r_{\alpha'}$ and $r'$, but this is a contradiction because $(y_\alpha,b_\alpha)\subset f^{-1}[(y,+\infty)]$.  
\end{proof}

\begin{corollary}
	Let \( f:[0,1]\rightarrow\Rea_\kappa \) be \( \kappa \)-continuous, and let \( x \in [0,1] \), \( \seq{r_\alpha}_{\alpha<\kappa} \) and \( \seq{r'_\alpha}_{\alpha <\kappa} \) be respectively increasing and decreasing sequences in \( [0,1] \) such that for all \( \alpha < \kappa \) we have \( f(r_\alpha) < x \) and \( f(r'_\alpha) > x \).
	Then there exists \( y \in [0,1] \) such that \( f(y) = x \) and \( \intset{r_\alpha}{\alpha < \kappa} < \extset{y} < \intset{r'_\alpha}{\alpha < \kappa} \).
	\label{kappa-continuous_preserves_value_in_between}
\end{corollary}
\begin{proof}
	Construct a sequence \( \seq{x_\alpha}_{\alpha < \gamma} \) for some \( \gamma \leq \kappa \) as follows.
	First let \( \delta_0 = 1 \).
	Having constructed \( \seq{x_\beta}_{\beta < \alpha} \) for some even \( \alpha < \kappa \), by Lemma~\ref{Theo:continuity} there exists \( x_\alpha \in [0,1] \) such that \( f(x_\alpha) = x \) and \( \intset{r_\beta}{\beta < \sup_{\nu < \alpha} \delta_\nu } < \extset{x_\alpha} < \intset{r'_\beta}{\beta < \sup_{\nu < \alpha} \delta_\nu } \).
	If \( \intset{r_\beta}{\beta < \kappa} < \extset{x_\alpha} < \intset{r'_\beta}{\beta < \kappa} \), then we are done and \( \gamma = \alpha \).
	Otherwise there exists \( \beta < \kappa \) such that \( r_\beta > x \) or \( r'_\beta < x \), so we let \( x_{\alpha+1} = r_\beta \) or \( x_{\alpha+1} = r'_\beta \) accordingly, and let \( \delta_\alpha = \beta+1 \).
	If the construction goes on for \( \kappa \) steps, then \( \seq{x_\alpha}_{\alpha<\kappa} \) is as in Lemma~\ref{kappa_continuous_cannot_flip_kappa_times}, a contradiction.
	Hence the construction ends at some stage \( \gamma < \kappa \), and therefore \( \intset{r_\beta}{\beta < \kappa} < \extset{x_\gamma} < \intset{r'_\beta}{\beta < \kappa} \).
\end{proof}

\begin{theorem}\label{Theo:main}
\begin{enumerate}
\item 
	If there exists an effective enumeration of a dense subset of $\Rea_\kappa$, then $\mathrm{IVT}_\kappa \leq_{\mathrm{W}} \mathrm{B}^{\kappa}_\mathrm{I}$. \label{IVT_B}
\item 
	We have $\mathrm{B}^{\kappa}_\mathrm{I} \leq_{\mathrm{W}} \mathrm{IVT}_\kappa$. \label{IVT_BCont}
\item 
	We have $\mathrm{IVT}_\kappa \leq^\mathrm{t}_{\mathrm{W}} \mathrm{B}^{\kappa}_\mathrm{I}$, and therefore $\mathrm{IVT}_\kappa \equiv^\mathrm{t}_{\mathrm{W}} \mathrm{B}^{\kappa}_\mathrm{I}$.\label{IVTEQ}
\end{enumerate}
\end{theorem}
\begin{proof}
	For item \ref{IVT_B}, let the $\kappa$-continuous function $f:[0,1]\rightarrow\Rea_\kappa$ be given, $\mathbb{D}$ be a dense subset of $\Rea_\kappa$ and
 \( \seq{d_\gamma}_{\gamma<\kappa} \) be an effective enumeration of \( [0,1] \cap \mathbb{D} \).
	Without loss of generality we can assume $f(0)<0$ and $f(1)>0$, and start setting $r_0=0$ and $r'_0=1$.
	Now assume that for $0<\alpha<\kappa$ we have already defined an increasing sequence $\seq{r_\beta}_{\beta<\alpha}$ and a decreasing sequence $\seq{r'_\beta}_{\beta<\alpha}$ of elements of \( [0,1] \cap \mathbb{D} \) with \( \intset{r_\beta}{\beta<\alpha} < \intset{r'_\beta}{\beta<\alpha} \) and \( \intset{f(r_\beta)}{\beta<\alpha} < \extset{0} < \intset{f(r'_\beta)}{\beta<\alpha} \).
	By Lemma \ref{Theo:continuity} there is still a root of \( f \) between the two sequences.
	Note that, since \( \Rea_\kappa \) is an \( \eta_\kappa \)-set and again by applying Lemma \ref{Theo:continuity}, there exist \( r_L,r_R \in \mathbb{D} \) such that \( \intset{r_\beta}{\beta<\alpha} < \extset{r_L} < \extset{r_R} < \intset{r'_\beta}{\beta<\alpha} \) and \( f(r_L) < 0 \), \( f(r_R) > 0 \).
	Therefore, by searching in the sequence \( \seq{d_\gamma}_{\gamma < \kappa} \) and running the corresponding algorithms in parallel, we can find such a pair \( r_L,r_R \) in fewer than \( \kappa \) computation steps.
	Let \( \beta, \gamma, \delta \) be such that \( \Goedel(\beta,\Goedel(\gamma,\delta)) = \alpha \), where \( \Goedel \) is the G\"odel pairing function, which has a computable inverse by Lemma~\ref{Goedel_pairing}.
	If \( r_L < d_\gamma < d_\delta < r_R \), \( f(d_\gamma) < 0 \), and \( f(d_\delta) > 0 \), where the last two comparisons are decided in fewer than \( \beta \) steps of computation, then let \( r_\alpha = d_\gamma \) and \( r'_\alpha = d_\delta \); otherwise let \( r_\alpha = r_L \) and \( r'_\alpha = r_R \).

	By Corollary~\ref{kappa-continuous_preserves_value_in_between} we have that there exists \( x \in [0,1] \) such that \( \intset{r_\alpha}{\alpha < \kappa} < \extset{x} < \intset{r'_\alpha}{\alpha < \kappa} \).
	It remains to be proved that \( f(x) = 0 \) for any such \( x \).
	Suppose not, say \( f(x) > 0 \) for some such \( x \).
	Then also \( f(y) > 0 \) for some \( y \in \mathbb{D} \) such that  \( \intset{r_\alpha}{\alpha < \kappa} < \extset{y} < \intset{r'_\alpha}{\alpha < \kappa} \).
	Now let \( \beta,\gamma,\delta < \kappa \) be such that  \( d_\gamma = y \),  \( d_\delta = r_\nu \) for some \( \nu \) such that \( \extset{y -_\mathrm{s} r_\nu} < \intset{r'_\alpha -_\mathrm{s} r_\beta}{\alpha,\beta < \kappa} \) and \( f(y) < 0 \), \( f(r_\nu) > 0 \) are decided in fewer than \( \beta \) computation steps.
	Then at stage \( \alpha = \Goedel(\beta,\Goedel(\gamma,\delta)) \) of the computation we define a pair \( r_\alpha, r'_\alpha \) such that \( r'_\alpha -_\mathrm{s} r_\alpha \leq y -_\mathrm{s} r_\nu \), a contradiction. This ends the proof of \ref{IVT_B}. 
	
Item \ref{IVT_BCont} can be proved by a straightforward generalisation of the proof of \cite[Theorem 6.2]{Brattka}, and the proof of item \ref{IVTEQ} is the same as that of item \ref{IVT_B} without the requirement that the enumeration $\langle d_\gamma \rangle_{\gamma < \kappa}$ of the dense subset of $[0,1] \cap \mathbb{D}$ be effective.

\end{proof}

Note that the antecedent of item \ref{IVT_B} of Theorem \ref{Theo:main} is satisfied, e.g., in the constructible universe $\mathbf{L}$. 
We leave for future work the task of investigating the set-theoretic properties of that condition more deeply.

\subsubsection*{Acknowledgments}

The authors would like to thank the Isaac Newton Institute for Mathematical Sciences
for the hospitality during the research programme Mathematical, Foundational
and Computational Aspects of the Higher Infinite.
The research benefited from the Royal Society International Exchange Grant \emph{Infinite games in logic and Weihrauch degrees}. 
The second author was also supported by the Capes Science Without Borders grant number 9625/13-5. 
The authors are grateful to Benedikt L\"{o}we and Arno Pauly for the many fruitful discussions and to the Institute for Logic, Language and Computation for the hospitality offered to the first author.
Finally, the authors wish to thank the three anonymous referees for the helpful comments which have improved the paper.

\end{document}